\documentclass[11pt]{amsart}
\usepackage{amssymb}
\newtheorem{theorem}{{\bf Theorem}}[section]
\newtheorem{lemma}{{\bf Lemma}}[section]
\newtheorem{prop}{{\bf Proposition}}[section]

\begin{document}
\title{On the conjugacy problem for cyclic extensions of free groups}
\author[Bardakov]{Valerij Bardakov}
\address{Sobolev Institute of Mathematics, Novosibirsk 630090, Russia}
\email{bardakov@math.nsc.ru}
\author[Bokut]{Leonid Bokut}
\address{Sobolev Institute of Mathematics, Novosibirsk 630090, Russia}
\email{bokut@math.nsc.ru}
\author[Vesnin]{Andrei Vesnin}
\address{Sobolev Institute of Mathematics, Novosibirsk 630090,
Russia} \email{vesnin@math.nsc.ru}
%
%
\thanks{Authors were supported by RFBR (grant~02--01--01118).}
\date{\today}

\begin{abstract}
We study the conjugacy problem in cyclic extensions of free
groups. It is shown that the conjugacy problem is solvable in
split extensions of finitely generated free groups by virtually
inner auto\-mor\-phisms. An algorithm for construction of the
unique representative (the conju\-ga\-cy normal form) for each
conjugacy class is given.
\end{abstract}

\maketitle


\section{Introduction}

One of fundamental problems in the combinatorial group theory is
the conjugacy problem formulated by M.~Dehn (1912) together with
two other problems: the word problem and the isomorphism problem
\cite[Ch.~2, \S~1]{LS}. M.~Dehn proved that the conjugacy problem
is solvable in fundamental groups of closed orientable surfaces.
Later his method have been extended to the class of groups with
small cancellation \cite[Ch.~5]{LS}.

Let us recall basic results about solvability of the conjugacy
problem in some classes of groups.

In \cite{Gromov} M.~Gromov introduced a class of groups which are
now referred to as word hyperbolic groups. Among examples of word
hyperbolic groups are finite groups, free groups, small
cancellation groups satisfying a metric small cancellation
condition $C'(\lambda)$ with $0 < \lambda \leq 1/6$. In
particular, the fundamental group of an oriented surface of genus
$g>1$ is hyperbolic. It is known \cite{BH} that the conjugacy
problem in word hyperbolic groups is solvable. Also it is solvable
in such generalizations of word hyperbolic groups as relative
hyperbolic groups and semi-hyperbolic groups. At the same time,
word hyperbolic groups are contained in the class of bi-automatic
groups, which are contained in the class of automatic groups
itself. The class of automatic groups belong to the class of
combable group. It is known that the conjugacy problem is solvable
in bi-automatic groups. Recently, M.R.~Bridson \cite{Bridson}
demonstrated that there exist combable groups in which the
conjugacy prob\-lem is unsolvable. The question about solvability
of the conjugacy problem in automatic groups is still open.

We will interested in groups $F_n(t) = F_n \leftthreetimes \langle
t\rangle$ which are semi-direct pro\-ducts of a free group $F_n$
and a cyclic group $\langle t \rangle$, where conjugation by $t$
induces an automorphism $\varphi \in \text{Aut} (F_n)$. In
particular, if $t$ is of infinite order, then $F_n(t)$ is the
mapping torus of $F_n$ corresponding to the automorphism $\varphi$
and is denoted by
$$
F_n (\varphi) \, = \, \langle x_1, x_2, \ldots, x_n, t \, | \,
t^{-1} x_i t = \varphi (x_i), \quad i=1,2,\ldots, n \rangle.
$$
If $t$ is of finite order (which is divided by order of $\varphi$)
then there exists a homomorphism of $F_n (\varphi)$ onto $F_n(t)$.

The study of groups $F_n (\varphi)$ is also motivated by their
relation with fundamental groups of closed 3-manifolds fibering
over a circle (see \cite{FH}).

In \cite{GS} J.M.~Gersten and  J.R.~Stallings asked about word
hyperbolicity of $F_n (\varphi)$. It was shown by  M.~Bestvina and
M.~Feighn \cite{BF} and by P.~Brinkmann \cite{Br}
that $F_n (\varphi)$ is word hyperbolic if and
only if $\varphi$ having no nontrivial periodic conjugacy classes
(i.e. there no non-trivial element $f$ of $F_n$ and non-zero
integer $k$ such that $\varphi^k(f)$ is conjugated to $f$ in
$F_n$). The solvability of the conjugacy problem for such groups
follows from the solvability of the conjugacy problem in word
hyperbolic groups.

We remark that some one-relator groups as well as some Artin
groups can be presented as cyclic extensions of free groups (see
discussions in Section~2). The conjugacy problem in one-relator
groups in still open. At the same time, it is solvable in
one-relator groups with torsion. This result was announced by
B.~Newman \cite{Newman}. S.~Pride \cite{Pride} proved this fact
for the case when the defining relation is of the form $r^n$ for
$n>2$. Another proof of this fact was given by V.N.~Bezverhnii
\cite{B}. L.~Larsen \cite{Larsen1} proved that the conjugacy
problem is solvable in one-relator groups with non-trivial center.
For some classes of one-relator groups the conjugacy problem was
solved by G.A.~Gurevich \cite{Gurevich1, Gurevich2} and
A.A.~Fridman \cite{Fridman}.

It is known that the conjugacy problem is solvable in braid groups
and knot groups \cite{Garside, Sela, Weinbaum}. The solvability of
the conjugacy problem in link groups seems still open.

Recall that the conjugacy problem is solvable in the Novikov group
$A_{p_1,p_2}$ \cite{Novikov} if and only if the word problem is
solvable in the corresponding Post system $P(A_{p_1, p_2})$
\cite{Bokut68}, \cite[Ch.~7]{BK}.

Surveys on the conjugacy problem in various classes of groups can
be found in \cite{Hurwitz84} and \cite{NRR}. The problem on
solvability of the conjugacy problem in groups $F_n(\varphi)$ was
formulated also by I.~Kapovich \cite[Problem~6.2]{Kapovich}.

It is known that the solvability of the conjugacy problem does not
preserve under finite extensions \cite{GK}. Moreover, D.~Collins
and C.~Miller \cite{CM} constructed a group $G$ containing a
subgroup $H$ of index two such that the conjugacy problem is
solvable in $H$ but not solvable in $G$. At the same paper they
constructed a group with solvable conjugacy problem which contains
a subgroup of index two with non-solvable conjugacy problem.

One of possible approaches to solve the conjugacy problem in
groups $F_n(t)$ is investigation of the property to be conjugacy
separable.

A group $G$ is said to be {\it conjugacy separable} if for any
pair of non-conjugated elements of $G$ there exists a homomorphism
of $G$ to a finite group such that images of these elements are
also non-conjugated. It was shown by A.I.~Malcev \cite{Malcev1},
if a finitely generated group is conjugacy separable then the
conjugacy problem in this group is solvable. Thus, the following
problem aries naturally:

\underline{Problem.} {\em Is a group $F_n(t)$ conjugacy
separable?}

It was shown by J.L.~Dyer \cite{Dyer1} that finite extensions of
free groups are conjugacy separable. Also, he proved in
\cite{Dyer2} that if $G$ is an extension of a free group by a
cyclic group and center of $G$ is nontrivial, then $G$ is
conjugacy separable. Moreover, if $G$ is one-relator group with
non-trivial center, then $G$ is conjugacy separable.

In the present paper we consider groups $G = F_n \leftthreetimes
\langle  t \rangle$ which are semi-direct product of a free group
$F_n$ and a cyclic group $\langle t \rangle$ such that the
conjugation by $t$ induces automorphism $\varphi \in \text{Aut}
(F_n)$ such that $\varphi^m \in \text{Inn} (F_n)$, i.~e.
$\varphi^m$ is inner automorphism of $F_n$ for some positive
integer $m$. Such automorphism $\varphi$ will be referred to as
{\it virtually inner}. In Section~2 we will demonstrate that many
one-relator groups and, in particular, two-generator Artin groups
can be obtained as split extension of free groups by virtually
inner automorphisms (see Section~2). In Section 3 and Section 4 we
will show that the conjugacy problem is solvable in $G$. For each
element of $G$ we will construct unique conjugacy normal from. If
$t$ is of finite order, then $G$ is word hyperbolic, and moreover,
in virtue of the above referred result of J.L.~Dyer, it is
conjugacy separable. Therefore, the conjugacy problem is solvable
in $G$. But our approach gives the more effective solving
algorithm than the general solving algorithm for word hyperbolic
groups. If $t$ is of infinite order, then $G$ is not word
hyperbolic because it contains a subgroup isomorphic to
$\mathbb{Z} \times \mathbb{Z}$ (see [19]). In Section~5 we will
show that the conjugacy problem is solvable for the split
extension $F_{\infty} \, \ltimes \, \mathbb Z$ of the countably
generated free group $F_{\infty}$ by the automorphism $\varphi$
shifting its generators.

\section{One-relator groups and 2-generated Artin groups}

In this section we will show that some one-relator groups and, in
par\-ti\-cular, 2-generator Artin groups can be obtained as split
extensions of free groups by virtually inner automorphisms.

It was shown in \cite{Appel, Bez, Bez-Dobr} that the conjugacy
problem is solvable in Artin and Coxeter groups of large type.
Also, it was shown in \cite{BS} that it is solvable for Artin
groups of finite type, i.e. such that corresponding Coxeter groups
are finite. The question about solvability of the conjugacy
problem in arbitrary Artin group seems still open.

Consider a group
$$
G \, = \, \langle  t, a \, | \, r = 1 \rangle
$$
with two generators and one defining relation. Let us assume that
the word $r$ is cyclically reduced and contains both generators
$a$ and $t$.

Applying, if necessary, Lemma~11.8 from \cite[ch.~5]{LS} we can
assume that the exponent sum of $t$ in $r$ is equal to zero. We
will use the Magnus--Moldavanskii method to represent $G$ as
HNN-extension.  Let us define new generators $a_i = t^{-i} a t^i$,
$i \in \mathbb{Z}$ and let $r'$ be the presentation of the word
$r$ written in these generators. We write $a_i \in r'$ if $a_i$ or
$a_i^{-1}$ is a subword of $r'$, and denote $\mu = \min \{ i \, |
\,  a_i \in r'\}$ and $\nu = \max \{ i \, | \, a_i \in r'\}$. Let
us define
$$
H  \, = \, \langle  a_i, \, \mu \leq i \leq \nu \, | \, r' = 1
\rangle ,
$$
$$
A \, = \, \langle a_i, \, \mu \leq i < \nu \rangle,
$$
$$
B \, = \, \langle a_i, \, \mu < i \leq \nu  \rangle.
$$
Then we can write $G = \langle H, t \, | \, t^{-1} A t = B,
\varphi \rangle$, where isomorphism $\varphi$ acts as the
following
$$
\varphi(a_i) = a_{i+1},~~~i = \mu, \mu + 1,\ldots, \nu - 1.
$$
If $a_{\nu}$ appears in $r'$ only once (with degree $+1$ or $-1$),
then using $r'$ we can express $a_{\nu}$ via other generators and
eliminate it from the generating system for $H$. Thus $H$ is a
free generated group which coincides with $A$. In this case
$$
G \, = \, \langle H, t \, | \,  t^{-1} H t = B, \varphi \rangle
$$
is said to be an {\it ascending HNN-extension} of $H$. If,
moreover, $B = H$, then $\varphi$ is an automorphism of $H$, and
$G = F_n (\varphi)$, where $n = \nu - \mu -1$, is a cyclic
extension of a free group.

Analogously, if $a_{\mu}$ can be expressed from $r'$ via other
generators, we get that $H$ coincides with $B$.

\smallskip

As a noticeable  example, let us consider 2-generated Artin group:
\begin{equation*}
\mathcal A(m) = \langle x, y  \, | \, w_m (x, y) = w_m (y,x )
\rangle ,
\end{equation*}
where $m \geq 3$  is integer and
$$
w_m (u,v) \, = \,
\left\{
\begin{array}{ll}
(u v)^n , &  \mbox{if} \quad m = 2n, \\
(u v)^n \, u , &   \mbox{if} \quad m = 2n+1 .
\end{array}
\right.
$$
We recall that $\mathcal A(2n+1)$ is the fundamental group of the
$(2n+1,2)$-torus knot complement in the 3-sphere and $\mathcal
A(2n)$ is the fundamental group of the $(2n, 2)$-torus 2-component
link complement in the 3-sphere. These knots and links arise as
closures of 2-strand braids.

Changing generators of $\mathcal A(m)$ in the same way as in
\cite{BV} (where the Gr\"obner--Shirshov bases for these groups
were constructed) we will get the following result.

\begin{prop} \label{P2.1}
Any 2-generated Artin group is a split extension of a free group
of finite rank and a cyclic group generated by a virtually inner
automorphism.
\end{prop}

\begin{proof}
Let us consider the group $\mathcal A(2n) = \langle x, y \, | \,
(xy)^n = (yx)^n \rangle $, $n \geq 2$, Denoting $t=x$ and $y_i =
t^i y t^{-i}$, where $i=0,1,\ldots, n-1$, we get the following
presentation:
\begin{equation*}
\mathcal A(2n)  = F_n (\varphi)  = \langle y_0, y_1, \ldots,
y_{n-1}, t \, | \, t^{-1} y_i t  =  \varphi(y_i) \, , \quad i=0,
\ldots, n-1 \rangle ,
\end{equation*}
where $F_n$ is the free group generated by $y_0, \ldots,
y_{n-1}$ and $\varphi$ is defined by:
\begin{eqnarray*}
\varphi(y_0) & = & y_0 y_1 \cdots y_{n-2} y_{n-1} y_{n-2}^{-1}
\cdots y_1^{-1} y_0^{-1}, \\
\varphi(y_i) & = & y_{i-1}, \qquad \qquad i=1, 2, \ldots, n-1.
\end{eqnarray*}
It is easy to check that the element $\Delta  = y_0 y_1 \cdots
y_{n-1}$ is such that $\varphi (\Delta) = \Delta$  (see also
\cite{BV}) and the automorpihsm $\varphi^n$ acts as the following:
$$
\varphi^n (y_i) \, = \, \Delta \, y_i \, \Delta^{-1}, \qquad i =
0, 1, \ldots , n-1,
$$
so $\varphi$ is an virtually inner automorphism of $F_n$.

Let us consider the group $\mathcal A(2n+1) = \langle x, y \, | \,
(x y)^n x = (y x)^n y \rangle $. Denoting $t=x$, $z = y x^{-1} = y
t^{-1}$, and $z_i = t^i z t^{-i}$, where $i=0,1,\ldots, 2n-1$, we
get the following presentation
\begin{equation*}
\mathcal A(2n+1) = F_{2n} (\psi) = \langle z_{0},  \ldots,
z_{2n-1}, t \, | \, t^{-1} z_i t = \psi(z_i) , \quad i=0, \ldots,
2n-1 \rangle ,
\end{equation*}
where $F_{2n}$ is the free group generated by $z_0, z_1, \ldots,
z_{2n-1}$ and $\psi$ is defined by:
\begin{eqnarray*}
\psi (z_0) & = & z_0 z_2 \cdots z_{2n-2}  z_{2n-1}^{-1}
z_{2n-3}^{-1} \cdots z_3^{-1} z_1^{-1} , \\
\qquad \psi (z_i) & = & z_{i-1}, \qquad \qquad i=1, \ldots, 2n-1.
\end{eqnarray*}
It is easy to check that the element
$$
\Sigma \, = \, z_0 z_2 \cdots z_{2n-2} \, (z_0 z_1 \cdots z_{2n-1})^{-1}
\, z_1 z_3 \cdots z_{2n-1}
$$
is such that $\psi (\Sigma) = \Sigma$ (see also \cite{BV}) and the
automorphism $\psi^{2(2n+1)}$ acts as the following:
$$
\psi^{2(2n+1)} (z_i) \, = \,  \Sigma \, z_i \, \Sigma^{-1},
\qquad i = 0, 1, \ldots , 2n- 1,
$$
so, $\psi$ is an virtually inner automorphism of $F_n$.
\end{proof}

We remark that the conjugacy problem in 2-generator Artin groups
is solvable since these groups are Artin groups of finite type.

\section{Finitely generated free groups and virtually inner automorphisms}

Let $F_n = \langle x_1, x_2, \ldots, x_n \rangle$ be the free
group of rank $n\ge 2$ with words from the alphabet ${\mathbb X} =
\{ x_1^{\pm 1}, \ldots, x_n^{\pm 1}\}$. In some cases, we will
need to distinguish words which represent the same element of the
group. We will write $U = V$ if two words (or two elements of the
group) are equal as elements of the group, and $U \equiv V$ if two
words are equal graphically. Denote by $|V|$ the length of a word
$V$ in the alphabet $\mathbb X$. A word $V$ is said to be {\it
reduced} if it contains no part $x x^{-1}$, $x \in \mathbb X$. A
reduced word $V$ defines a non-identity element if and only if
$|V| \ge 1$. A reduced word obtained by reducing of an original
word will be referred to as its {\it reduction}. By $|| V ||$ we
will denote the length of the reduction of a word $V$. Further, a
reduced word $V=x_{i_1}^{\varepsilon_1} x_{i_2}^{\varepsilon_2}
\cdots x_{i_n}^{\varepsilon_n},$ where $\varepsilon_i = \pm 1,$
$i=1, \ldots, n$, is said to be {\it cyclically reduced} if $i_1
\neq i_n$ or if $i_1 = i_n$ then $\varepsilon_1 \neq -
\varepsilon_n$. Clearly, every element of a free group is
conjugated to an element given by a cyclically reduced word
referred to as its {\it cyclic reduction}.

We will use standard notations $\text{\rm Aut} (F_n)$ and
$\text{\rm Inn} (F_n)$ for the group of automorphisms and the
group of inner automorphisms of $F_n$, respectively. An
automor\-phism $\varphi \in \text{\rm Aut} (F_n)$ is said to be
{\it virtually inner} if $\varphi^m \in \text{\rm Inn} (F_n) $ for
some positive integer $m$.

In particular, any automorphism $\varphi$ of finite order is
virtually inner.

Denote by $F_n(t) =  F_n \, \ltimes \, \langle t \rangle$ the
semi-direct product, where the ge\-ne\-ra\-tor $t$ of the cyclic
group $\langle t \rangle$ is such that the conjugation by $t$
induces an automorphism $\varphi$, i.e. $t^{-1} f t = \varphi (f)$
for any $f \in F_n$. Remark that $F_n \triangleleft F_n(t)$.

We will show that the following property holds.

\begin{theorem}\label{T3.1}
If $\varphi$ is a virtually inner automorphism of a  free group
$F_n$, $n \geq 2$, then the conjugacy problem is solvable for
$F_n(t)$.
\end{theorem}

To prove this result, we will find an unique representative for
each conju\-ga\-cy class.

Any element $v \in F_n(t)$  can be presented in the form $t^{\ell}
V$ where $\ell \in \mathbb Z$ if $t$ has infinite order, and $0
\leq \ell <  |t|$ if $t$ has finite order $|t|$; $V$ is a word of
the alphabet $\mathbb X$ (we will also say that $V$ is a $\mathbb
X$-{\it part} of $v$). Moreover, if $V$ is reduced then such a
presentation of $v$ is unique.

Let $\varphi \in \text{\rm Aut} (F_n)$ be a virtually inner
automorphism such that $\psi = \varphi^m$ is an inner automorphism
of $F_n$. Without loss of generality we can assume that $m$ is
taken the smallest positive integer having such a property. There
exists a reduced word $\Delta \in F_n$ such that
$$
\varphi^m (f) \, = \, \Delta^{-1} f \Delta
$$
for any $f \in F_n$.
It is easy to check (see also
\cite[Lemma~2]{BBV}) that following properties hold.

\begin{lemma}   \label{L3.1}
(1) Automorphisms $\varphi$ and $\psi$ commute.\\ (2) If $k = m q
+ r$, where $q \in \mathbb Z$ and $0\leq r \leq m-1$, then for any
word $U$ of the alphabet $\mathbb X$ we have
$$
\varphi^k (U) \, = \, \varphi^r (\Delta^{-q} \, U \, \Delta^q).
$$\\
(3) $\varphi (\Delta) = \Delta$.
\end{lemma}

Let us define a linear order ``$<$''  on the set of irreducible
words in the alphabet $\mathbb X$. Assume that elements of
$\mathbb X$ are ordered in the following way:
$$
x_1 \, < \, x_1^{-1} \, < \, x_2 \, < \, x_2^{-1} \, < \, \ldots
\, < \, x_n \, < \, x_n^{-1} .
$$
We write $U < V$ if $|U| < |V|$ or if $|U| = |V|$ and the word $U$
is less than the word $V$ in respect to the lexicographical order
corresponding to the above defined linear order on $\mathbb X$.

A reduced word $V \in F_n$ is said to be {\it $\Delta$-reduced} if
$|V | \leq || \Delta^{-k} V \Delta^k ||$ for all $k \in \mathbb
Z$. Obviously, if $V$ is cyclically reduced, then the length of
any word conjugated to $V$ is not less than the length of $V$, so
$V$ is $\Delta$-reduced.

\begin{lemma}\label{L3.2}
Suppose that $\Delta$ is cyclically reduced. A reduced word $V \in
F_n$ is $\Delta$-{\it reduced} if $|V| \leq ||
\Delta^{-\varepsilon} V \Delta^\varepsilon ||$ for $\varepsilon =
\pm 1$.
\end{lemma}

\begin{proof} See \cite[Lemma~3]{BBV}.
\end{proof}

If $\Delta$ is cyclically reduced, Lemma~\ref{L3.2} gives the
finite algorithm to find for a given reduced word $V$ a
$\Delta$-reduced word $V_{\Delta}$ conjugated to $V$ by some power
of $\Delta$. Indeed, it is enough to repeat conjugations of $V$ by
$\Delta^{\varepsilon}$, $\varepsilon = \pm 1$, few times. If the
length of the obtained word is less than the length of the
previous word, we will conjugate again. If not, then the obtained
word is a $\Delta$-reduced word $V_{\Delta}$ conjugated to $V$.
Such a construction of a $\Delta$-reduced word $V_{\Delta}$
conjugated to $V$ by some power of $\Delta$ will be referred to as
a $\Delta$-{\it reduction}.

We remark that if $\Delta$ is not cyclically reduced, then the
analog of Lemma~\ref{L3.2} does not hold. It is clear from the
following example.

\smallskip

{\bf Example.} Let $U, W, \Sigma \in F_n$ be nonempty reduced
words such that for an integer $|k|>1$ words $\Delta \equiv U^{-1}
W^{-1} \Sigma W U$ and $V \equiv U^{-1} W^{-1} \Sigma^k W U^2$ are
reduced. If $k>1$, we get
\begin{eqnarray*}
\Delta^{-1} V \Delta & = & U^{-1} W^{-1} \Sigma^{k-1} W U W^{-1}
\Sigma W U, \cr \Delta^{-2} V \Delta^2 & = & U^{-1} W^{-1}
\Sigma^{k-2} W U W^{-1} \Sigma^2 W U, \cr & \cdots & \cr
\Delta^{-(k-1)} V \Delta^{k-1} & = & U^{-1} W^{-1} \Sigma W U
W^{-1} \Sigma^{k-1} W U, \cr \Delta^{-k} V \Delta^k & = & W^{-1}
\Sigma^k W U. \end{eqnarray*}
It is easy to see that
$$
|V| < ||\Delta^{-1} V \Delta|| = ||\Delta^{-2} V \Delta^2 || =
\ldots = ||\Delta^{-(k-1)} V \Delta^{k-1}||,
$$
but $||\Delta^{-k} V \Delta^k || < |V|$.

If $k < -1$, similar example can be obtained.

\smallskip

\begin{lemma} \label{L3.3}
Suppose $\Delta$ is not cyclically reduced. Let $V$ be a reduced
word such that $|V| \leq ||\Delta^{-\varepsilon} V
\Delta^{\varepsilon} ||$, $\varepsilon = \pm 1$. Then one of the following cases
holds:\\
(1) $V$ is $\Delta$-reduced.  \\
(2) $\Delta \equiv U_1^{-1} U_2^{-1} \Delta_{11}^{-1} \Delta_2
\Delta_{11} U_2 U_1$ and $V \equiv U_1^{-1} U_2^{-1} W U_2 U_1$,
for some reduced $\Delta_{11}, \Delta_{2}, U_1, U_2$ and
cyclically reduced $W$ for which there exist integers $k$,
$|k|>1$, and
$m \geq 0$ and reduced $\Phi$, $W_0$ such that either\\
a) $U_2^{-1} \Delta_{11}^{-1} \Delta_2^{-k} \Delta_{11} \equiv
\Phi W^{-m}$ and $W \equiv \Phi^{-1} W_0$,\\
or\\
b) $\Delta_{11}^{-1} \Delta_2^{k} \Delta_{11} U_2 \equiv W^{-m}
\Phi$ and $W \equiv W_0 \Phi^{-1}$.\\
In addition, $\Phi$ does not end by $W^{-1}$, and $W_0$ is either
nonempty or $W_0 \equiv 1$ with $U_1^{-1} \Phi^{-1} U_1$ be
reduced.

In case (2) $V_{\Delta} = \Delta^{-k} V \Delta^k$ is
$\Delta$-reduced word conjugated to $V$. Moreover, case (2)
describe all possible cases when $|V|\leq ||\Delta^{-\varepsilon}
V \Delta^{\varepsilon}||$, $\varepsilon = \pm 1$, but
$||\Delta^{-k} V \Delta^k|| < |V|$ for some $|k|>1$.
\end{lemma}

\begin{proof} See \cite[Lemma~4]{BBV}.
\end{proof}

If $\Delta$ is not cyclically reduced, Lemma~\ref{L3.3} gives the
finite algorithm to find for a given reduced word $V$ a
$\Delta$-reduced word $V_{\Delta}$ conjugated to $V$ by some power
of $\Delta$. If $\Delta$ and $V$ are of the form represented in
case (2) of Lemma~\ref{L3.3}, then we define $V_{\Delta} =
\Delta^{-k} V \Delta^k$. In this case we say that $V_{\Delta}$ is
a {\it $\Delta$-reduction} of $V$. If $\Delta$ and $V$ are others,
then we follow the same steps as described after Lemma~\ref{L3.2}.

Remark that $V_{\Delta}$ is not uniquely determined by $V$ (see
\cite[Proposition~2]{BBV}).

\section{Constructing of conjugated normal form}

Now we will study $\Delta$-reduced elements of the group $F_n(t)$.
A word $v = t^{\ell} V \in F_n(t)$, $\ell \in \mathbb Z$, $V \in
F_n$, is said to be {\it $\Delta$-reduced} if $V$ is
$\Delta$-reduced.

In virtue of Lemma  \ref{L3.1}(3), the conjugation of $v =
t^{\ell} V$ by any power $\Delta^m$ means the conjugation of $V$
by $\Delta^m$, i.e.
\begin{eqnarray*}
\Delta^{-m} v \Delta^m & = & \Delta^{-m} t^{\ell} V \Delta^m \quad
= \quad t^\ell (t^{-\ell} \Delta^{-m} t^\ell ) V \Delta^m \cr & =
& t^\ell \varphi^{\ell} (\Delta^{-m}) V \Delta^m \quad = \quad
t^\ell \Delta^{-m} V \Delta^m .
\end{eqnarray*}

For any integer $k$ we will denote by $V_{[k]}$ a reduction of
$t^{-k} V t^k = \varphi^k (V)$. Remark that
$t^{-m}V_{[l]}t^m=V_{[l+m]}$ and that $V_{[0]}=V$ if $V$ is a
reduced word.

If length of a reduced word $V$ in the alphabet $\mathbb X$ is
bigger than $1$ then it can be presented as a product $V \equiv V'
V''$ of two nonempty reduced words. The word $V'$ will be referred
to as an {\it initial part} of $V$ and $V''$ will be referred as a
{\it final part} of $V$. Denote by $\mathcal I(V)$ the set of all
initial parts of $V$ and by $\mathcal F(V)$ the set of all final
parts of $V$.

Consider $v \in F_n(t)$ and fix its presentation $v \equiv
t^{\ell} V$, $\ell \in \mathbb Z$, where $V$ is reduced word in
the alphabet $\mathbb X$. Conjugation by $t^\ell$ induces an
automorphism $\psi = \varphi^\ell$ of the free group $F_n$.

For a reduced word $V \equiv V' V''$ in the alphabet $\mathbb X$
with $V' \in \mathcal I(V)$, $V'' \in \mathcal F (V)$ we say that
a word $V''_{[\ell]} V'$ is a {\it cyclic $\psi$-shift of a final
part of $V$} and a word $V'' V'_{[-\ell]}$ is a {\it cyclic
$\psi$-shift of an initial part of $V$}. Remark that these
elements can be obtained by a conjugation of $v$. Indeed,
conjugating $v$ by $(V'')^{-1}$ we will get
$$
V'' v (V'')^{-1} \, = \,
t^\ell t^{-\ell} V'' t^\ell V' \, = \,
t^\ell V''_{[\ell]} V',
$$
and conjugating $v$ by $(V')_{[-\ell]}$ we will get
$$
(V'_{[-\ell]})^{-1} v V'_{[\ell]} \, = \,
t^\ell t^{-\ell} (V'_{[-\ell]})^{-1} t^\ell V' V'' V'_{[-\ell]} \, = \,
t^\ell V'' V'_{[-\ell]} .
$$
If $|V'| = 1$, i.e. $V' = x_i^{\varepsilon}$, $\varepsilon = \pm
1$, then the corresponding cyclic $\psi$-shift will be referred to
as a {\it cyclic $\psi$-shift of the initial letter}. If  $|V''| =
1$ then the corresponding cyclic $\psi$-shift will be referred to
as a {\it cyclic $\psi$-shift of the final letter}.

A reduced word $V$ will be referred to as a {\it cyclically
$\psi$-reduced} if there doesn't exist a cyclic $\psi$-shift of
its final part and there doesn't exist a cyclic $\psi$-shift of
its initial part which decreases length of $V$. Obviously,
applying cyclic $\psi$-shifts of final (or initial) parts to a
given word $V$ we will get a cyclically $\psi$-reduced word
conjugated to $V$.

Now let us construct conjugating normal form for an element $v
 = t^\ell V$. Without loss of generality, we can assume that
$V$ is cyclically $\psi$-reduced and $\Delta$-reduced. For a given
$V$ let us construct a set of words $\mathcal V_{\Delta}$ which
contains $V$ and such elements which are conjugated to $V$ by
elements of the group $\langle \Delta \rangle$ and are
$\Delta$-reduced. All of them have the same length as $V$. By
Lemma~\ref{L3.2} and Lemma~\ref{L3.3} the set $\mathcal
V_{\Delta}$ is finite. Applying to elements of $\mathcal
V_{\Delta}$ cyclic $\psi$-shifts of all initial parts and all
final parts, we will construct the set $(\mathcal
V_{\Delta})_{\psi}$. Applying to obtained words $\Delta$-reducing,
we will get a set $((\mathcal V_{\Delta})_{\psi})_{\Delta}$ each
element of which is $\Delta$-reduced. Consider the subset of
$((\mathcal V_{\Delta})_{\psi})_{\Delta}$ consisting of cyclically
$\psi$-reduced words. Multiplying each of obtained words on left
by $t^{\ell}$ we will get a set of words from $F_n(t)$. Let us
denote the obtained set by $D_0(v)$.

Remark that in the free group the set of words obtained from a
given word $V$ by cyclic shifts is finite. But the set of words
obtained from a word $v \in F_n(t)$ by cyclic $\psi$-shifts can be
infinite if $t$ has infinite order. Indeed, it is clear from
following relations
\begin{eqnarray}
\psi^{k-1}(V) \, \psi^{k-2}(V) \, \cdots \, \psi(V) \, V \, v \, V^{-1} \,
\psi(V^{-1}) \, \cdots \, \psi^{k-2}(V^{-1}) \, \psi^{k-1}(V^{-1}) \nonumber \\
\, = \, t^{\ell} \, \psi^{k}(V),
\qquad \mbox{for} \qquad  k>0; \nonumber
\end{eqnarray}
\begin{eqnarray}
\psi^{k}(V^{-1}) \, \psi^{k+1}(V^{-1}) \, \cdots \, \psi^{-1}(V^{-1}) \, v \,
\psi^{-1}(V) \, \cdots \, \psi^{k+1}(V) \, \psi^{k}(V) \nonumber \\
\, = \, t^{\ell} \, \psi^{k}(V), \qquad \mbox{for} \qquad k<0, \nonumber
\end{eqnarray}
that the element $v\equiv t^{\ell}V$ is conjugated to
$t^{\ell}\psi^k(V)$ for any integer $k$. Moreover, the conjugation
can be done by an element of the free group $F_n$.

For each integer $k$ we define a set $D_k(v) = D_0 (t^{\ell}
\psi^k(V))$, and $D(v) = \cup_{k \in \mathbb Z} D_k(v)$. Let us
verify that $D(v)$ is finite.

\begin{lemma} \label{L4.1}
For integer $\ell$ and $m$ as above denote $d = \gcd (m, \ell)$. Then
$$
D(v) \, = \, \cup_{k \in \{ 0, d, \ldots, m-d\}} D_0 (t^{\ell}\varphi^k(V)) .
$$
\end{lemma}

\begin{proof}
For any integer $k$ we have $\psi^k (V) = \varphi^{\ell k} (V)$.
Let $r$, $0\leq r < m$ be such that $\ell k = m q + r$, $q\in
\mathbb Z$. By Lemma~\ref{L3.1},
$$
\varphi^{\ell k} (V) \, = \, \Delta^{-q} \, \varphi^r (V) \, \Delta^q .
$$
Let $\ell = \ell_1 d$ and $m = m_1 d$ for some integer $\ell_1$
and $m_1$. If $k$ runs over the set $\mathbb Z$ then $r$ runs over
the set $\{ 0, d, 2d, \ldots, m-d \}$. Let us show that $D_k(v) =
D_0 (t^{\ell}\varphi^r(V))$, that will give the statement. Indeed, by the definition,
 $D_k (v) = D_0 (t^{\ell}
\Delta^{-q} \varphi^r (V) \Delta^q)$. Denote $U = V_{[r]} =
\varphi^r (V)$ and consider elements from
$D_0(t^{\ell}\varphi^r(V))$. By the definition, $D_0
(t^{\ell}\varphi^r(V))$ consists of words whose $\mathbb X$-parts
are the $\psi$-shifts $\psi(U_2) U_1$ or $U_2 \psi^{-1} (U_1)$,
where $U \equiv U_1 U_2$, to which $\Delta$-reducing and
$\psi$-reducing are applied. To construct $D_k(v)$ we must pass
from a word $\Delta^{-q} \, \varphi^r (V) \, \Delta^q$ to a
$\Delta$-reduced word, which, in a general case, can be different
from $U = \varphi^r (V)$. But, according to the definition of $D_0
(v)$, the set of $\Delta$-reduced words constructed from
$\Delta^{-q} \, \varphi^r (V) \, \Delta^q$ coincides with the set
of $\Delta$-reduced words constructed from $\varphi^r(V)$. So,
corresponding sets of all cyclic $\psi$-shifts of initial parts
and of finals parts also coincide. Therefore, $D_k (v) = D_0
(t^{\ell}\varphi^r(V))$.
\end{proof}

\begin{lemma}\label{L4.2}
The set $D(v)$ has the following properties: \\
(1) If $u \in D(v)$ then $D(u) = D(v)$; \\
(2) Let elements $v \equiv t^{\ell} V$ and $w \equiv t^{\ell} W$ be
conjugated by an element of the group $H = \langle F_n, t^m
\rangle$. Suppose that words $V$ and $W$ are cyclically
$\psi$-reduced and $\Delta$-reduced. Then $D(v) = D(w)$.
\end{lemma}

\begin{proof}
Since $D(v)$ consists of elements conjugated by elements of $H$,
item (2), obviously, implies (1). Let us prove (2). Let $u \equiv
t^{ms} U$ be a conjugating element, where $s$ is some integer and
$U$ is a reduced word. By the assumption, $w = u^{-1} v u$, i.e.
we have the following equality in the free group $F_n$:
\begin{equation}
W \, = \, U^{-1}_{[\ell]} \, \Delta^{-s}  \, V \, \Delta^{s} \, U
. \label{eq1}
\end{equation}
Since $t^{\ell} \Delta^{-s} =\Delta^{-s}  t^{\ell}$, denoting $U' =
\Delta^s U$ we get $U^{-1}_{[\ell]} \Delta^{-s} =
(U'_{[\ell]})^{-1}$, i.e. $W= U^{-1}_{[\ell]} \Delta^{-s} V
\Delta^{s} U = (U'_{[\ell]})^{-1} V U'$. Therefore, without loss
of generality, we can assume that in (\ref{eq1}) the exponent $s$
is equal to zero. Since $W$ is cyclically $\psi$-reduced, the
product $U^{-1}_{[\ell]} V U$ contains cancellations. Moreover,
these cancellations are either in the product $U^{-1}_{[\ell]} V$
or in the product $V U$, but not in the both.

\underline{Case 1.} Suppose that $V$ is not cancelling wholly.

\underline{Case 1(a).} Suppose that there are cancellations in the
product $U^{-1}_{[\ell]} V$. Then $U^{-1}_{[\ell]}$ must be
cancelled wholly, i.e. $V \equiv U_{[\ell]} V_1$. Then
$$
W \, = \, U^{-1}_{[\ell]} V U \, \equiv \, U^{-1}_{[\ell]} (U_{[\ell]} V_1) U
\, =  \, V_1 U
$$
and $w \equiv t^{\ell} V_1 U$ is obtained from $v \equiv t^{\ell} U_{[\ell]} V_1$
by the shift of the initial part $U_{[\ell]}$. Indeed, according to the above
described procedure, we need to conjugate $v$ by $U = U_{[0]}$:
$$
U^{-1}_{[0]} v U_{[0]} \, = \,
U^{-1}_{[0]} (t^{\ell} U_{[\ell]} V_1) U_{[0]} \, = \,
t^{\ell} U^{-1}_{[\ell]} U_{[\ell]} V_1 U_{[0]} \, = \,
t^{\ell} V_1 U = w .
$$
Therefore, $D(v) = D(w)$.

\underline{Case 1(b).} Suppose that there are cancellations in the
product $VU$. Then $U$ must be cancelled wholly, i.e. $V = V_1
U^{-1}$. Then
$$
W \, = \, U^{-1}_{[\ell]} V U \, \equiv \, U^{-1}_{[\ell]} ( V_1 U^{-1} ) U
\, = \, U^{-1}_{[\ell]} V_1
$$
and $w \equiv t^{\ell} U^{-1}_{[\ell]} V_1$ is obtained from $v \equiv
t^{\ell} V_1 U^{-1}$ by the shift of the final part $U^{-1}$.
Indeed, conjugating $v$ by $U$ we get
$$
U^{-1} v U \, \equiv \, U^{-1} (t^{\ell} V_1 U^{-1}) U \,= \,
U^{-1}_{[\ell]} V_1 .
$$
Therefore, $D(v) = D(w)$.

\underline{Case 2.} Suppose that $V$ is cancelling wholly.

\underline{Case 2(a).} Suppose that $V$ is cancelling wholly in
the product $V U$ and after that there are cancellations of
letters of the remaining part of $U$ with letters of
$U^{-1}_{[\ell]}$. Then we can represent $U \equiv V^{-1} U_1$,
therefore $U^{-1}_{[\ell]} = U^{-1}_{1 \, [\ell]} V_{[\ell]}$ and
\begin{equation}
U^{-1}_{[\ell]} V U \, = \,  U^{-1}_{1 \, [\ell]} V_{[\ell]} V V^{-1} U_1
\, = \, U^{-1}_{1 \, [\ell]} V_{[\ell]} U_1 .
\label{eq3}
\end{equation}
If $U_1$ is cancelling wholly with $V_{[\ell]}$ then $V_{[\ell]} \equiv V_2
U_1^{-1}$ and
$$
U^{-1}_{1 \, [\ell]} V_{[\ell]} U_1 \, \equiv \,
U^{-1}_{1 \, [\ell]} V_2 U^{-1}_1 U_1 \, = \,
U^{-1}_{1 \, [\ell]} V_2 \, =  \, W .
$$

If word $V_{[\ell]}$ in the product $V_{[\ell]} U_1$ is cancelling
wholly, then we will use induction by length of $U$.

Remark that $W$ arises in the process of the construction of $D_1
(v)$. Indeed,
$$
D_1(v) = D_0 (t^{\ell} V_{[\ell]}) = D_0 (t^{\ell} V_2 U_1^{-1})
$$
and $W$ is obtained from $V_2 U_1^{-1}$ by the cyclic $\psi$-shift
of the final part.

\underline{Case 2(b).} Suppose that $V$ is cancelling wholly in
the product $U^{-1}_{[\ell]} V$ and after that there are
cancellations of letters of the remaining part of
$U^{-1}_{[\ell]}$ with letters of $U$. Then we can represent
$U^{-1}_{[\ell]} \equiv U_1^{-1} V^{-1}$, therefore $U = \psi^{-1}
(V) \psi^{-1} (U_1)$ and
\begin{equation}
U^{-1}_{[\ell]} V U = U^{-1}_1 V^{-1} V V_{[-\ell]} U_{1 \, [-\ell]} \, = \,
U^{-1}_1 V_{[-\ell]} U_{1 \, [-\ell]} . \label{eq2}
\end{equation}
If $U_1^{-1}$ is cancelling wholly with $V_{[-\ell]}$ then
$V_{[-\ell]} = U_1 V_2$ and $U_1^{-1} V_{[-\ell]} U_{1 \, [-\ell]}
= U_1^{-1} (U_1 V_2) U_{1 \, [-\ell]} = V_2 U_{1 \, [-\ell]}$,
i.e. $V = U_{1 \, [\ell]} V_{2 \, [\ell]}$ and $W \equiv V_2 U_{1
[-\ell]}$. Comparing these words we see that $W$ belongs to
$D_{-1}(v)$. Indeed, by the definition,
$$
D_{-1} (v) \, = \, D_{0} (t^{\ell} \psi^{-1} (V)) \, = \, D_0 (t^{\ell} U_1
V_2) .
$$
Applying to $U_1 V_2$ the cyclic $\psi$-shift of the initial part
$U_1$, we will get $W = V_2 U_{1 \, [-\ell]}$. The case when
$V_{[-\ell]}$ is cancelling wholly with $U^{-1}_1$ in (\ref{eq2})
can be considered similar to the above.
\end{proof}

Let us define
$$
\overline{D} (v) \, = \, \cup_{0\le k < m} \, D(t^{-k} v t^k) .
$$
Obviously, this set is finite. Let $\overline{v} = t^{\ell} V_0$
be an element with the smallest $\mathbb X$-part among all
elements of $\overline{D}(v)$. Such $\overline{v}$ is said to be
the {\it conjugacy normal form} for $v$. By the construction, $v$
and $\overline{v}$ are conjugated in $F_n(t)$. To show uniqueness
of $\overline{v}$ we will use the following statement.

\begin{lemma}\label{L4.3}
Let $w \equiv t^{\ell} W \in F_n(t)$ be cyclically $\psi$-reduced
and $\Delta$-reduced. If $w$ is conjugated to $v$ in $F_n(t)$ then
$\overline{D} (w) = \overline{D} (v)$.
\end{lemma}

\begin{proof}
Suppose that $v$ and $w$ are conjugated in $F_n(t)$ by $u \equiv
t^k U$, where $k \in \mathbb Z$ and $U$ is a reduced word in the
alphabet $\mathbb X$. Let $r$,
 $0 \leq r \leq m-1$, be such that $k = m q + r$ for some integer $q$.
 By Lemma~\ref{L3.1} we have
\begin{eqnarray*}
u^{-1} v u & = & U^{-1} t^{-k} (t^{\ell} V) t^k U \quad = \quad
U^{-1} t^{\ell} (\Delta^{-q} t^{-r} V t^r \Delta^q ) U \cr & = &
U^{-1} t^{-mq} t^{\ell} (t^{-r} V t^r) t^{mq} U = U^{-1} t^{-mq}
(t^{-r} v t^r) t^{mq} U ,
\end{eqnarray*} i.e. $w$ is conjugated to $t^{-r} v t^r$ by an element from the group
$\langle F_n, t^m\rangle$. By the construction,  $D(t^{-r} v t^r)
\subseteq \overline{D} (v)$, and by Lemma~\ref{L4.2},
$$
D(t^{-r} v t^r) = D (U^{-1} t^{-mq} (t^{-r} v t^r) t^{mq} U).
$$
Therefore, $\overline{D}(w) = \overline{D}(v)$.
\end{proof}

Now we are able to complete the proof of Theorem~\ref{T3.1}

\begin{proof}
Let $v = t^{\ell} V$ be an element of $F_n(t)$. Using, if
necessary, conjugation by elements of $F_n$, we can assume that
$V$ is cyclically $\psi$-reduced and $\Delta$-reduced. Let us
construct a set $\overline{D} (v)$ as in Section~4. From this set
we choose words of minimal length. Then, from such words, choose
the conjugacy normal form $\overline{v}$ for $v$, that is the word
whose $\mathbb X$-part is minimal in respect to the above ordering
on $F_n$.

For a pair of given words $u = t^k U$ and $v = t^{\ell} V$ the
conjugacy problem is solving as the following. If $k \neq \ell$
then $u$ and $v$ are not conjugated in $F_n(t)$. If $k = \ell$
then let us construct conjugacy normal forms $\overline{u}$ and
$\overline{v}$. By Lemma~\ref{L4.3} words $u$ and $v$ are
conjugated in $F_n(t)$ if and only if $\overline{u} =
\overline{v}$.
\end{proof}

We remark that if $t$ is of finite order in $F_n(t)$, then this
group is almost free, so, it is word hyperbolic. This implies the
solvability of the conjugacy problem in this group. But our
approach gives the more effective solving algorithm than the
general solving algorithm for word hyperbolic groups.

\smallskip

The following question naturally arises in the context of the above obtained result.

\underline{Problem:} {\em Let $G=F \, \leftthreetimes M$, where
$M\subseteq \mbox{\rm Aut}(F_n)$ is such that the image of $M$ in
the group $\mbox{\rm Out}(F_n) = \mbox{\rm Aut} (F_n) / \mbox{\rm
Inn}(F_n)$ under the natural homomorphism is finite group. Does
the conjugacy problem solvable in $G$?}

\smallskip

Let us show that the answer on this question is affirmative if $M$ is finite group.

\begin{prop} \label{P4.1}
Let $G = F_n \, \leftthreetimes M$, where $M$ is a finite subgroup
of $\text{\rm Aut } (F_n)$. Then the conjugacy problem is solvable
for $G$.
\end{prop}

\begin{proof}
Let us suppose that $M$ has $k$ elements:
$$
M=\{ \alpha_0=e, \alpha_1,\ldots, \alpha_k \} .
$$
We order these elements in the following way:
$$
\alpha_0 < \alpha_1 < \ldots < \alpha_k.
$$
Consider elements $v=\alpha_i V$ and $u=\alpha_j U$ from $G$,
where $V$ and $U$ are reduced words from $F_n$. Obviously, if
elements $\alpha_i$ and $\alpha_j$ are not conjugated in $M$, then
elements $v$ and $u$ are not conjugated in $G$. So, we can suppose
that $\alpha_i$ and $\alpha_j$ are conjugated in $M$. Using, if
necessary, a conjugation, we can assume that $v=\varphi V$ and
$u=\varphi U$, where $\varphi \in M$, and $\varphi$  is taken the
smallest representative of the class of conjugated elements that
contains $\alpha_i$ and $\alpha_j$. Further, as in the proof of
Theorem~\ref{T3.1},  we construct sets $D_0(v)$ and $D_0(u)$. In
addition, $\Delta \equiv 1$. Define
$$
D(v)=\bigcup_{i=0}^{k}D_0(v^{\alpha_i}), \qquad
D(u)=\bigcup_{i=0}^{k}D_0(u^{\alpha_i}).
$$
For each of these sets we choose a word that has the smallest
$\mathbb X$-part in respect to the above defined ordering. Such a
word will be the conjugated normal form of the corresponding word.
Therefore, if obtained words coincide, then elements $u$ and $v$
are conjugated in $G$. If they are different, then elements $u$
and $v$ are not conjugated in $G$.
\end{proof}

\section{Countably generated free group and shifting
automorphism}

Let $F_{\infty} = \langle \ldots , x_{-1}, x_0, x_1, x_2, \ldots
\rangle$ be the free group with countably infinite number of
generators with words from the alphabet $\mathbb X = \{ \ldots,
x_{-1}^{\pm 1}, x_0^{\pm 1}, \\ x_1^{\pm 1}, x_2^{\pm 1}, \ldots
\}$. Consider an automorphism $\varphi : F_{\infty} \to
F_{\infty}$ acting by shifting generators: $\varphi (x_i) =
x_{i-1}$, where $i \in \mathbb Z$. We will call $\varphi$ the {\it
shifting automorphism}. Denote by $F_{\infty} (\varphi)  =
F_{\infty} \ltimes \langle t \rangle$ the split extension, where
$t$ is the generator of the cyclic group $\langle t \rangle$ such
that the conjugation by $t$ induces the automorphism $\varphi$,
i.e. $t^{-1} x_i t = x_{i-1}$ for $i \in \mathbb Z$. Remark that
$F_{\infty} \triangleleft F_{\infty} (\varphi)$.

Defining relations for $F_{\infty} (\varphi)$ can be written as
following:
$$
\begin{array}{ll}
x_i t \, = \, t x_{i-1} , & \qquad x_i t^{-1} = t^{-1} x_{i+1}  \\
x^{-1}_i t \, = \, t x^{-1}_{i-1} , & \qquad x^{-1}_i t^{-1} \, =
\, t^{-1} x_{i+1}^{-1},
\end{array}
$$
where $i \in \mathbb Z$. Using the same arguments as in the proof
of Lemma~2.1 from \cite{BV}, we remark that these relations
together with trivial relations form the Gr\"obner--Shirshov basis
for $F_{\infty} (\varphi)$.

In the present section we will show that the following property holds.

\begin{theorem}\label{T5.1}
If $\phi$ is the shift automorphism of the free group $F_{\infty}$
then the conjugacy problem is solvable for $F_{\infty} (\varphi)$.
\end{theorem}

Any element of $F_{\infty}(\varphi)$ is uniquely presented in the
form $t^m V$, where $m \in \mathbb Z$ and $V \in F_{\infty}$ is a
reduced word. To prove the decidability of the conjugacy problem
for the group $F_{\infty}(\varphi)$ we will show that for each
class of conjugated elements we can choose the unique
representative (the conjugacy normal form), and that two elements
from $F_{\infty}(\varphi)$ are conjugated if and only if the
representatives of corresponding classes coincide.

For any integer $m$ let $\varphi^m$ be an automorphism of the
group $F_{\infty}$ defined by
$$
\varphi^m (x_i) \, = \, t^{-m} x_i t^m = x_{i-m} .
$$
By $V_{[m]}$ we denote a free reduction of the word $\varphi^m
(V)$. A word $W$ is said to be {\it $\varphi^m$-conjugated} to a
word $V$ by a word $U^{-1}$ if
$$
W \, \equiv \, \varphi^m (U) V U^{-1} \, = \, U_{[m]} V U^{-1}.
$$
Obviously, two words $t^m V$ and $t^m W$ from
$F_{\infty}(\varphi)$ are conjugated by an element of $F_{\infty}$
if and only if corresponding words $V$ and $W$ from $F_{\infty}$
are $\varphi^m$--conjugated. A reduced word  $V$ is said to be
{\it cyclically $\varphi^m$-reduced} if it can not be presented in
the form
$$
V \equiv \psi^m(U) \, V_0  \, U^{-1} ,
$$
for non-trivial $U \in F_{\infty}$.

\begin{lemma}\label{L5.1}
Each element $v \in F_{\infty}(\varphi)$, presented by a word $t^m
V$, is conjugated to a word $t^m V_0$, where $V_0$ is cyclically
$\varphi^m$-reduced.
\end{lemma}

\begin{proof}
If $V$ is not cyclically $\varphi^m$-reduced, then $v$ can be
presented in the form
$$
v \equiv t^m \, \varphi^m (U) \, V_0 \, U^{-1} ,
$$
where $V_0$ is cyclically $\varphi^m$-reduced Conjugating this
word by $U$ we get
$$
U^{-1} v U  \equiv  U^{-1} \, ( t^m \varphi^m (U) V_0 U^{-1} ) \,
U \, = \, U^{-1} t^m t^{-m} U t^m  V_0 \, = \, t^m V_0 .
$$
\end{proof}

A word $v \equiv t^m \, V \in F_{\infty}(\varphi)$ is said to be
{\it cyclically $\varphi^m$-reduced} if the word $V \in
F_{\infty}$ is cyclically $\varphi^m$-reduced.

Suppose that a word $V\in F_{\infty}$ has some $x_k$ as the final
letter, i.e. $V \equiv U x_k$. Let us define the {\it cyclic
$\varphi^m$-shift of the final letter\,}, denoted by $\tau_m$, as
following:
$$
\tau_m \, : \, V \, \equiv \, U x_k \, \mapsto \, \varphi^m (x_k)
\, V \, x_k^{-1} \, \equiv \, \varphi^m (x_k) \, ( U x_k ) \,
x_k^{-1} \, = \, x_{k-m} U.
$$
Obviously, if $v = t^m V$ then $v' = t^m \tau_m (V)$ is conjugated to $v$
by $x_j^{-1}$:
$$
x_j v x_j^{-1} \, = \, x_j ( t^m U x_j) x_j^{-1} \, = \,
t^m x_j t^{-m} t^m U \, = \, t^m x_{j-m} U \, = \, v' .
$$

Let us define a linear order "$<$" on the set of reduced words of
the alphabet $\mathbb X$. Assume that elements of $\mathbb X$ are
ordered in the following way:
$$
\ldots < x_{-1} < x_{-1}^{-1} < x_0 < x_0^{-1} < x_1 < x_1^{-1} <
x_2 < x_2^{-1} <  \ldots
$$
We write $U < V$ if $|U|<|V|$ or if $|U|=|V|$ and the word $U$ is
less than the word $V$ in respect to the lexicographical order
corresponding to the above defined linear order on $\mathbb X$.

Remark the following obvious property:

\begin{lemma}\label{L5.2}
Let words $U_1, U_2, \ldots, U_p$ representing elements of $F_{\infty}$
be ordered in the following way:
$$
U_1 \, < \, U_2 \, < \, \ldots \, < \, U_p ,
$$
and $m$ be an integer. Then the following ordering holds
$$
\varphi^m (U_1) \, < \, \varphi^m (U_2) \, < \, \ldots \, < \,
\varphi^m (U_p) ,
$$
that can be also written as
$$
U_{1 \, [m]} \, < \, U_{2 \, [m]} \, < \, \ldots \, < \, U_{p \, [m]} .
$$
\end{lemma}

\medskip

For a $\varphi^m$-reduced word $v \equiv t^m V \in
F_{\infty}(\varphi)$, with $|V| = n$,  define a set
$$
D (v) \, = \, \{ t^m \, \tau^i_m (V) \quad | \quad i=0,1,\ldots, n-1 \} .
$$
Let us choose a word from $D(v)$, say $t^m V_0$, such that
its $\mathbb X$-part $V_0$ is smallest in respect to above defined
lexicographical order.
Suppose that $V_0$ starts from a letter $x_k^{\varepsilon}$,
$\varepsilon = \pm 1$. Then the word $\varphi^k (V_0)$ starts from
the letter $x_0^{\varepsilon}$. For such chosen $V_0$ and $k$, the
word $\overline{v} = t^m \varphi^k (V_0)$ will be referred to as the {\it
conjugacy normal form} for $v$. Since the automorphism $\varphi^k$
acts by the conjugation and the $\varphi^m$--conjugation can be
realized by the conjugation too, any word is conjugated to its
 conjugacy normal form.

\begin{lemma}\label{L5.3}
If words $v, w \in F_{\infty}(\varphi)$ are conjugated, then their
conjugacy normal forms coincide.
\end{lemma}

\begin{proof}
Suppose that we have found the conjugacy normal form for a given
element $v = t^{k} V$. By Lemma~\ref{L5.1}, we can assume, up to a
conjugation, that $v$ is taken to be cyclically
$\varphi^k$-reduced. Let us consider an element conjugated to $v$
and demonstrate that their conjugacy normal forms coincide.
Indeed, let $u=t^{\ell} U$ and
$$
w = u^{-1} v u \equiv U^{-1} t^{-\ell} t^{k} V t^{ \ell} U =
t^{k} \, U^{-1}_{[k]} \, V_{[\ell]} U .
$$
Consider possible cases of cancellations in the word
$U^{-1}_{[k]} \, V_{[\ell]} \, U$.

\underline{Case 1.} Suppose that there are no cancellations. Then
$$
U \, w \, U^{-1} \, = \,
U \, t^k \, U^{-1}_{[k]} \, V_{[\ell]} \, = \, t^k \, V_{[\ell]} .
$$
Thus, for any $m$, the cyclic $\varphi^m$-shifts of $U w U^{-1}$
can be obtained from the cyclic $\varphi^m$-shift of $v$ using
conjugation by $t^{-\ell}$. Therefore, a conjugate normal form of
$w$ coincides with a conjugate normal form of $v$.

\smallskip

Let us assume that below $w\equiv t^k W$ is cyclically
$\varphi^k$-reduced. Using induction by length $|U|$ of  $U$ we will
show that $W$ can be obtained from $V$ by cyclic
$\varphi^k$--shift and conjugation by some degree of $t$. Then
there must be some cancellations in the word $U^{-1}_{[k]} \,
V_{[\ell]} \, U$.

\underline{Case 2.} Suppose that there are cancellations in the
product $V_{[\ell]} \, U$, i.e. $V_{[\ell]} \, \equiv \,
V'_{[\ell]} \, V''_{[\ell]}$ and $U \, \equiv \, (V''_{[\ell]}
)^{-1} \, U'$. Then
$$
U^{-1}_{[k]} \, V_{[\ell]} \, U \, \equiv \,
(U')^{-1}_{[k]} \, V''_{[\ell+k]} \, ( V'_{[\ell]} \, V''_{[\ell]}
\,)
(V''_{[\ell]})^{-1} \, U' \, = \,
(U'_{[k]})^{-1} \, V''_{[\ell+k]} \, V'_{[\ell]} \, U' .
$$
Remark that there can not be cancellations in the product
$V''_{[\ell+k]} \, V'_{[\ell]}$ because, by the assumption,
$v=t^{k} V$ is cyclically $\varphi^k$-reduced. Remark that the
element $V''_{[\ell+k]} \, V'_{[\ell]}$ can be obtained from $v$
by a cyclic $\varphi^k$--shift and a conjugation by $t^{\ell}$.
Since $|U'| < |U|$, the induction assumption can be applied.
Therefore, its normal form coincides with a conjugated normal form
of $v$.

\underline{Case 3.} Suppose that there are cancellations in the
product $U^{-1}_{[k]} \, V_{[\ell]}$, i.e.
$$
V_{[\ell]} \, \equiv \, V'_{[\ell]} \, V''_{[\ell]} ,
\qquad
U^{-1}_{[k]} \, = \, (U')^{-1} \, (V'_{[\ell]})^{-1} ,
\quad
U_{[k]} \, = \, V'_{[\ell]} \, U' .
$$
Then
$$
U^{-1}_{[k]} \, V_{[\ell]} \, U \, \equiv \, (U')^{-1} \,
(V'_{[\ell]})^{-1} \, ( V'_{[\ell]} \, V''_{[\ell]} \,)
V'_{[\ell-k]} \, U'_{[-k]} \, = \, (U')^{-1} \, V''_{[\ell]} \,
V'_{[\ell-k]} \, U'_{[-k]} ,
$$
and there are no cancellations in the product $V_{[\ell]} \, U$.
We see that the element $V''_{[\ell]} \, V'_{[\ell-k]}$ can be
obtained from $V' \, V''$ by a cyclic $k$--shift, i.e.
$$
V'' \, v \, (V'')^{-1} \, = \, t^{k} \, V''_{[k]} \, V',
$$
and a conjugation by $t^{-\ell+k}$. Since $|U'| < |U|$, the
induction assumption can be applied, and we get the statement.
\end{proof}

As a consequence of the above considerations we get the statement
of Theorem~\ref{T5.1}.



\begin{thebibliography}{99}

\bibitem{AD}
S.I.~Adyan, V.G.~Durnev, {\em Decision problems for groups and
semigroups,} Uspekhi Mat. Nauk {\bf 55} (2000), no. 2(332), 3--94
(Russian); translated in Russ. Math. Surveys {\bf 55} (2000),
no.~2, 207--296.

\bibitem{Appel}
K.~Appel, {\em On Artin groups and Coxeter groups of large type,}
in: Contributions in group theory, 50--78, Contemp. Math. {\bf
33}, Amer. Math. Soc., Providence, RI, 1984.

\bibitem{BBV}
V.~Bardakov, L.~Bokut, A.~Vesnin, {\em Twisted conjugacy in free
groups and Makanin's question,} Seoul National University, 2003,
RIM-GARC Preprint Series 03--12, 19 pp., aviliable from
arXiv:math.GR/0401349.

\bibitem{BFGMRRS}
K.~Bencsath, B.~Fine, A.M.~Gaglione, A.~Myasnikov, F.~Roehl, G.~Rosenberger,
D.~Spellman,
{\em Aspects of the theory of free groups,}
in: Algorithmic problems in groups and semigroups (Lincoln, NE, 1998), 51--90,
Trends Math., Birkha\"user Boston, Boston, MA, 2000.

\bibitem{BF}
M.~Bestvina, M.~Feighn, {\em A combination theorem for negatively
curved groups,} J. Diff. Geom., {\bf 35} (1992), no.~1, 85--101;
ibid. {\bf 43} (1996), no.~4, 783--788.


\bibitem{B}
V.N.~Bezverkhnii, {\it Solution of the conjugacy problem for words
in some classes of groups,} Algorithmic problems in the theory of
groups and semigroups, Tulsk. Gos. Ped. Inst., Tula,
1990, 103--152 (Russian).

\bibitem{Bez}
V.N.~Bezverkhnii, {\em Solution of the conjugacy problem for words
in Artin and Coxeter groups of large type,} Algorithmic problems in
the theory of groups and semigroups, Tulsk. Gos. Ped. Inst.,
Tula, 1986, 21--61 (Russian).

\bibitem{Bez-Dobr}
V.N.~Bezverkhnii, I.V.~Dobrynina, {\em Solution of the conjugacy
problem for words in Coxeter groups of large type}, Chebyshev Sbornik,
{\bf 4(1)} (2003), 10--33 (Russian).

\bibitem{Bokut68}
L.A.~Bokut, {\em Degrees of unsolvability of the conjugacy problem
for finitely presented groups,} Algebra i Logika {\bf 7} (1968),
no.~5, 4--70; ibid {\bf 7} (1968), no.~6, 4--52. (Russian)

\bibitem{BK}
L.A.~Bokut, G.P.~Kukin, {\em Algorithmic and combinatorial
algebra,} Mathematics and its Applications, {\bf 255}. Kluwer
Academic Publishers Group, Dordrecht, 1994.

\bibitem{BV}
L.~Bokut, A.~Vesnin, {\em Gr\"obner--Shirshov bases for some braid
groups,} 17~pp. to appear in Journal of Symbolic Computations,
(2004).

\bibitem{Bridson}
M.R.~Bridson, {\em The conjugacy and isomorphism problems for
combable groups,} Math. Ann. {\bf 327} (2003), no.~2, 305--314.

\bibitem{BH}
M.R.~Bridson, A.~Haefliger, {\em Metric spaces of non-positive
curvature,} Grundl. Math. Wiss., {\bf 319}, Springer-Verlag,
Berlin--Heidelberg, 1999.

\bibitem{BS}
E.~Brieskorn, K.~Saito, {\em Artin-Gruppen und Coxeter-Gruppen,}
Invent. Math. {\bf 17} (1972), 245--271.

\bibitem{Br}
P.~Brinkmann, {\em Hyperbolic automorphisms of free groups,} Geom.
Funct. Anal. {\bf 10} (2000), no.~5, 1071--1089.

\bibitem{CM}
D.~Collins, C.~Miller III, {\em The conjugasy problem and
subgroups of finite index,} Proc. London Math. Soc. (3) {\bf 34}
(1977), no.~3, 535--556.

\bibitem{Dehn}
M.~Dehn, {\em \"Uber unendichle diskontinuerliche Gruppen,} Math.
Ann. {\bf 71} (1912), 116--144.

\bibitem{Dyer1}
J.L.~Dyer, {\em Separating conjugates in free-by-finite groups,}
J. London Math. Soc. (2) {\bf 20} (1979), no.~2, 215--221.

\bibitem{Dyer2}
J.L.~Dyer, {\em Separating conjugates in amalgamated free products
and HNN extensions,} J. Austral. Math. Soc., Ser A {\bf 29}
(1980), no.~1, 35--51.

\bibitem{FH}
M.~Feighn, M.~Handel, {\em Mapping tori of free group
automorphisms are coherent,} Ann. of Math. (2) {\bf 149} (1999),
no.~3, 1061--1077.

\bibitem{Fridman}
A.A.~Fridman, {\em A solution of the conjugacy problem in a
certain class of groups,} in: Mathematical logic, theory of
algorithms and theory of sets (dedicated to P. S. Novikov on the
occasion of his seventieth birthday). Proc. Steklov Inst. Math.
{\bf 133} (1977), 233--242, 276.

\bibitem{Garside}
F.A.~Garside,
{\em The braid group and other groups,}
Quart. J. Math. Oxford Ser. (2) {\bf 20} (1969), 235--254.

\bibitem{GS}
S.M.~Gersten, J.R.~Stallings, {\em Irreducible outer automorphisms
of a free group,} Proc. Amer. Math. Soc., {\bf 111} (1991), no.~2,
309--314.

\bibitem{GK}
A.V.~Goryaga, A.S.~Kirkinskii,
{\em The decidability of the conjugacy problem cannot be transferred
to finite extensions of groups,}
Algebra i Logika {\bf 14(4)} (1975), 393--406. (Russian)

\bibitem{Gromov}
M.~Gromov, {\em Hyperblic groups,} in: Essays in Group Theory, Ed.
S.M.~Gersten, Math. Sci. Res. Inst. Publ. {\bf 8}, Springer,
New-York, 1987, 75--263.

\bibitem{Gurevich1}
G.A.~Gurevich, {\em On the conjugacy problem for groups with one
defining relation,} Sov. Math., Dokl. {\bf 13} (1972), 1436--1439.

\bibitem{Gurevich2}
G.A.~Gurevich, {\em On the conjugacy problem for groups with a
single defining relation,}  in: Mathematical logic, theory of
algorithms and theory of sets (dedicated to P. S. Novikov on the
occasion of his seventieth birthday). Proc. Steklov Inst. Math.
{\bf 133} (1977), 108--120.

\bibitem{Hurwitz84}
R.D.~Hurwitz,
{\em A survey of the conjugacy problem,}
Contributions to group theory,  278--298, Contemp. Math., {\bf 33},
Amer. Math. Soc., Providence, RI, 1984.

\bibitem{Kapovich}
I.~Kapovich, {\em Mapping tori of endomorphisms of free groups,}
Comm. in Algebra  {\bf 28} (2000), no.~6, 2895--2917.

\bibitem{KhS}
O.G.~Kharlampovich, M.V.~Sapir, {\em Algoritmic problems in
varieties,} Internat. J. Algebra Comput. {\bf 5} (1995), no.~4-5,
379--602.

\bibitem{Kour}
{\em The Kourovka notebook. Unsolved problems in group theory.}
Fifteenth augmented edition. Edited by V.D.~Mazurov and
E.I.~Khukhro. Russian Academy of Sciences Siberian Division,
Institute of Mathematics, Novosibirsk, 2002. 125 pp.

\bibitem{Larsen1}
L.~Larsen, {\em The conjugacy problem and cyclic HNN
constructions,} J. Austral. Math. Soc. Ser. A  {\bf 23} (1977),
no.~4, 385--401.

\bibitem{LS}
R.C.~Lyndon, P.E.~Schupp,
{\it Combinatorial group theory.}
Ergebnisse der Mathematik und ihrer Grenzgebiete, {\bf 89}.
Springer-Verlag, Berlin-New York, 1977.

\bibitem{Malcev1}
A.I.~Malcev, {\em On homomorphisms onto finite groups,} Uchen.
Zapiski Ivanovsk. ped. instituta {\bf 18} (1958), no.~5, 49--60
(also in ``Selected papers'', Vol. 1, Algebra, (1976), 450--461)

\bibitem{Malcev2}
A.I.~Malcev, {\em On isomorphi matrix representations of infinite
groups,} Mat. Sbornik {\bf 8} (1940), no.~3, 405--422.

\bibitem{Newman}
B.~Newman, {\em Some results on one-relator groups,} Bull. Amer.
Math. Soc. {\bf 74} (1968), 568--571.

\bibitem{NRR}
G.A.~Noskov, V.N.~Remeslennikov, V.A.~Roman'kov, {\em Infinite groups,}
Itogi Nauki Tekh., Ser. Algebra Topologiya Geom. {\bf 17} (1979), 65--157
(Russian).

\bibitem{Novikov}
P.S.~Novikov,
{\em Unsolvability of the conjugacy problem in the theory of groups,}
Izv. Akad. Nauk SSSR. Ser. Mat. {\bf 18} (1954), 485--524 (Russian).

\bibitem{Pride}
S.~Pride, {\it Small cancellation conditions satisfied by
one-relator groups,} Math. Z. {\bf 184} (1983), no.~2, 283--286.

\bibitem{Sela}
Z.~Sela, {\it The conjugcy problem for knot groups,} Topology {\bf
32} (1993), no.~2, 363--369.

\bibitem{Weinbaum}
C.M.~Weinbaum,
{\em The word and conjugacy problems for the knot group of any tame,
prime, alternating knot,}
Proc. Amer. Math. Soc. {\bf 30} (1971), 22--26.

\end{thebibliography}
\end{document}